\documentstyle[11pt,amstex,amssymb]{amsart}\textheight=20.1truecm\textwidth=13.7truecm
\begin{document} 

\ 

\vskip 30pt
\centerline{\Large\bf Matrix Versions of some Classical Inequalities}
\vskip 30pt
\centerline{Jean-Christophe Bourin}
\vskip 10pt
 \centerline{E-mail: bourinjc@@club-internet.fr}
 \vskip 5pt
 \centerline{8 rue Henri Durel, 78510 Triel, France}
 \vskip 10pt
 \centerline{\it Dedicated to Maslina Darus, with friendship and respect}
\vskip 20pt
\noindent
{\small {\bf Abstract.} 
\vskip 5pt
Some natural inequalities related to rearrangement in matrix products can also be regarded as extensions of 
classical inequalities for sequences or integrals. In particular, we show matrix versions of Chebyshev and
 Kantorovich type inequalities.
The matrix approach may also provide simplified proofs and new results for classical inequalities. For instance,
we show a link between Cassel's inequality and the basic rearrangement inequality for sequences of Hardy-Littlewood-Polya, and we state a reverse inequality to the Hardy-Littlewood-Polya inequality in which matrix technics are essential.


Mathematical subjects classification 2000:  15A60 47A30 47A63}

\vskip 25pt
{\large\bf Introduction}
\vskip 10pt
An important source of interesting inequalities in Matrix/operator theory is the study of rearrangements in a product.
An obvious, but useful, example is the operator norm inequality
\begin{equation}
\Vert AB\Vert_{\infty} \le \Vert BA\Vert_{\infty}
\end{equation}
whenever $AB$ is normal.
Here and in the sequel we use capital letters $A$, $B,\dots,\,Z$ to denote $n$-by-$n$ complex matrices, or operators on a finite dimensional Hilbert space ${\cal H}$; $I$ stands for the identity. When $A$ is positive semidefinite, resp.\ positive definite, we write $A\ge 0$, resp.\ $A>0$.
\vskip 10pt
In a series of papers, the author showed further rearrangement inequalities companion to (1). These inequalities may be considered as matrix versions of some classical inequalities for sequences or integrals.
The aim of this paper, which is mainly a survey, is  to emphasize the link between these classical inequalities and  matrix  rearrangement inequalities. The concerned classical inequalities are of Chebyshev and Kantorovich type. They can be stated for functions on a probability space, or equivalently for sequences. Let us present them for
 functions on an interval. In 1882, Chebyshev (see [11]) noted the following inequalities for bounded measurable functions $f$, $g$ on a real interval $\Omega$ endowed with a probability measure $\mu$: if $f$ and $g$ are both nondecreasing,
\begin{equation}
\int_{\Omega} f\,{\rm d}\mu \int_{\Omega} g\,{\rm d}\mu \le \int_{\Omega} fg\,{\rm d}\mu.
\end{equation}
Of course if $f$ and $-g$ are both nondecreasing then the reverse inequality holds. For measurable functions $f$ and $g$
with $p\ge f(t)\ge q$ and $r\ge g(t)\ge s$, Gruss showed in 1934 (see [21]) the following estimate for the difference in Chebyshev inequality,
\begin{equation}
\left|\int_{\Omega} f\,{\rm d}\mu \int_{\Omega} g\,{\rm d}\mu - \int_{\Omega} fg\,{\rm d} \mu\right| \le \frac{1}{4}(p-q)(r-s),
\end{equation}
in particular if $a\ge f(t)\ge b$,
\begin{equation}
\int_{\Omega} f^2\,{\rm d}\mu -\left(\int_{\Omega} f\,{\rm d}\mu\right)^2\le \,\frac{(a-b)^2}{4}.
\end{equation}
Such an inequality is called a Kantorovich type inequality. Indeed, when $b\ge0$, it is the additive version of
\begin{equation}
\int_{\Omega} f^2\,{\rm d}\mu \le \,\frac{(a+b)^2}{4ab}\left(\int_{\Omega} f\,{\rm d}\mu\right)^2
\end{equation}
which is equivalent to the original Kantorovich inequality stated in 1948 ([11]),
\begin{equation}
\int_{\Omega} f\,{\rm d}\mu \int_{\Omega} \frac{1}{f}\,{\rm d}\mu \le \,\frac{(a+b)^2}{4ab}.
\end{equation}

\vskip 15pt
The fact that some matrix inequalities can be regarded as generalizations of integral inequalities (2)-(6) takes roots in the observation that these integral inequalities have immediate operator reformulations.   By computing  inner products in an orthonormal  basis of eigenvectors for  $Z>0$ with extremal eigenvalues $a$ and $b$, the Kantorovich inequality  (6) may be rephrased as follows: for all norm one vectors $h$,
$$
\langle h, Zh\rangle \langle h, Z^{-1}h\rangle \le \,\frac{(a+b)^2}{4ab}.
$$
Similarly taking  square roots (5) can be rewritten 
\begin{equation}
\Vert Zh\Vert \le  \,\frac{a+b}{2\sqrt{ab}}\, \langle h, Zh\rangle .
\end{equation}
Such inner product inequalities are not only natural of their own right, they also motivate simple proofs via matrix techniques of the corresponding integral inequalities $-$ or their discrete analogous for sequences. For instance an extremely simple proof of (7) [5] (see also [6]) is reproduced here as Lemma 2.2.
A similar remark holds for a nice paper of Yamaziki [24] about the Specht reverse arithmetic-geometric mean inequality.

\vskip 10pt
Section 1  discusses an inequality for the Frobenius (or Hilbert-Schmidt) norm which is a matrix extension of  Chebyshev's inequality and which implies several classical inequalities, in particular von Neumann's Trace inequality. 
In the next section  we review several recent rearrangement inequalities for symmetric norms and eigenvalues which extend the Kantorovich inequalities. We also show a link between Cassel's inequality (a reverse Cauchy-Schwarz inequality) and a reverse inequality for the basic (Hardy-Littlewood-Polya) rearrangement inequality. 
Section 3 is concerned with matrix versions of generalized Kantorovich type inequalities.
\vskip 20 pt\noindent
{\large\bf 1. Matrix Chebyshev inequalities}
\vskip 10pt

The results of this section originate from [3] (see also [6]). Here we give simpler proofs based on some matrices with nonnegative entries associated to normal operators. We also show the connection with standard inequalities.

Say that a pair of Hermitians $(A,B)$ is monotone if there exists a third Hermitian $C$ and two nondecreasing functions $f$ and $g$ such that $A=f(C)$ and $B=g(C)$. If $g$ is nonincreasing  $(A,B)$ is antimonotone. We have the following result for the Frobenius (i.e., Hilbert-Schmidt) norm $\Vert\cdot\Vert_2$.

\vskip 10pt\noindent
{\bf Theorem 1.1.} {\it Let $A,\, B\ge 0$ and let $Z$ be normal.
If $(A,B)$ is monotone,
\begin{equation*}
\Vert AZB\Vert_2\le \Vert ZAB\Vert_2.
\end{equation*}
If $(A,B)$ is antimonotone,
\begin{equation*}
\Vert AZB\Vert_2\ge \Vert ZAB\Vert_2.
\end{equation*}
}

\vskip 5pt\noindent
Note that, for positive $A$ and $B$, $(A,B)$ is monotone iff so is $(A^{1/2},B^{1/2})$. Hence the first inequality  is equivalent to
\begin{equation*}
{\rm Tr}\, Z^*AZB \le {\rm Tr}\, Z^*ZAB.
\end{equation*}
Since for Hermitians $A$, $B$, the pair $(A,B)$ is monotone iff so is $(A+aI,B+bI)$ for any reals $a$, $b$, we then remark that Theorem 1.1 can be restated as trace inequalities involving Hermitian pairs:

\vskip 10pt\noindent
{\bf Theorem 1.2.} {\it Let $A$, $B$ be Hermitian and let $Z$ be normal.
If $(A,B)$ is monotone,
\begin{equation*}
{\rm Tr}\, Z^*AZB \le {\rm Tr}\, Z^*ZAB.
\end{equation*}
If $(A,B)$ is antimonotone,
\begin{equation*}
{\rm Tr}\, Z^*AZB \ge {\rm Tr}\, Z^*ZAB.
\end{equation*}
}

\noindent
 These inequalities have found an application to quantum information theory [13].

\vskip 10pt
Let us consider some classical facts related to these theorems.
The case of $Z$ unitary   entails von-Neumann's Trace inequality,
\begin{equation}
\left| {\rm Tr}\ XY\right| \le \sum \mu_j(X)\mu_j(Y)
\end{equation}
where $\mu_j(\cdot)$, $j=1,\dots$ denotes the singular values arranged in decreasing order and counted with their multiplicities. This important inequality is the key for standard proofs of the  H{\"o}lder inequality for Schatten $p$-norms. First if $X$ and $Y$ are both positive then, for some unitary $Z$, $ZXZ^*$ and $Y$ form a monotone pair and von Neumann's inequality follows from our theorems. For  general $X$ and $Y$,
 we consider  polar decompositions $X=U|X|$, $Y=V|Y|$ and we note that
\begin{align*}
\left| {\rm Tr}\, XY\right| &= \left| {\rm Tr}\, |X^*|^{1/2}U|X|^{1/2}V|Y|\right| \\
&= \left| {\rm Tr}\, |Y|^{1/2}|X^*|^{1/2}U|X|^{1/2}V|Y|^{1/2}\right| \\
&\le \{{\rm Tr}\, |Y||X^*|\}^{1/2} \{{\rm Tr}\, |X|V|Y|V^*\}^{1/2}
\end{align*}
by using the Cauchy-Schwarz inequality $\left| {\rm Tr}\, A^*B\right|\le \{{\rm Tr}\, A^*A\}^{1/2}
\{{\rm Tr}\, B^*B\}^{1/2}$ for all $A$, $B$. Hence the general case follows from the positive one.

Note that von Neumann's inequality is a matrix version of the classical Hardy-Littlewood rearrangement inequality: 
 Given real scalars $\{a_k\}_{k=1}^n$ and $\{b_k\}_{k=1}^n$, 
\begin{equation}
\sum_{k=1}^n a^{\uparrow}_k b^{\downarrow}_k \le   \sum_{k=1}^n a_k b_k \le \sum_{k=1}^n a^{\uparrow}_k b^{\uparrow}_k 
\end{equation}
where the exponent $\uparrow$ (resp.\ $\downarrow$) means the rearrangement in increasing (resp.\ decreasing) order.
These  inequalities also follow from Theorem 1.2 by letting $Z$ be a permutation matrix and $A$, $B$ be diagonal matrices.

Finally we remark that letting $Z$ be a rank one projection, $Z=h\otimes h$, then we get in the monotone case
\begin{equation}
\Vert Ah\Vert\,\Vert Bh\Vert\le \Vert ABh\Vert \quad{\rm and} \quad 
\langle h, Ah\rangle \langle h, Bh\rangle \le \langle h, ABh\rangle
\end{equation}
for all unit vectors $h$. the reverse inequalities hold for antimonotone pairs. These are just restatements of Chebyshev's inequality (2).
 
\vskip 10pt
The known proof of Theorems 1.1, 1.2  is quite intricated. The next lemma establishes a simple inequality for nonnegative matrices (i.e., with nonnegative entries) which entails the theorems.   Our motivation for searching a proof via nonnegative matrices was the following observation. If $Z=(z_{i,j})$ is a normal matrix, then $X=(x_{i,j})$ with $x_{i,j}=|z_{i,j}|^2$ is a  sum-symmetric matrix: for each index $j$,
\begin{equation*}
\sum_k x_{k,j} =\sum_k x_{j,k}.
\end{equation*}
Indeed, the normality of $Z$ entails $||Zh||^2=||Z^*h||^2$ for  vectors $h$, in particular  for  vectors of the canonical basis. 

For  a nonnegative matrix $X$, we define its row-column  ratio as the number
$$
{\rm rc}(X)=\max_{1\le i\le n}\frac{\sum_k x_{i,k}}{\sum_k x_{k,i}},
$$
whenever $X$ has at least one nonzero entry on each column. If not, we set  ${\rm rc}(X)=\lim_{r\to 0} {\rm rc}(X_r)$ where $X_r$ is the same matrix as $X$, except that the zero entries are replaced by $r>0$. It may happens that ${\rm rc}(X)=\infty$ and we adopt the convention $\infty\times 0=\infty$.
 
Given real column vectors,  $a=(a_1,\dots,a_n)^T$  and  $b=(b_1,\dots,b_n)^T$ we denote by   $a\!\cdot\!b$  the vector of the entrywise product of $a$ and $b$. We denote the sum of the components of the vector $a$ by $\sum a$.   We say that $a$ and $b$ form a monotone (resp. antimonotone) pair if, for all indexes $i$, $j$, 
\begin{equation*}
a_i<a_j\Rightarrow b_i\le  ({\rm resp.\ } \ge)b_j,
\end{equation*}
equivalently
\begin{equation*}
(a_i-a_j)(b_i-b_j)\ge ({\rm resp.\ } \le)0.
\end{equation*}
\vskip 10pt
We may then state results which compare $a\!\cdot\! X(b)$ and $X(a\!\cdot\! b)$\,:

\vskip 10pt\noindent
{\bf Lemma 1.3.} {\it Let $X$ be a nonnegative matrix and let $(a,b)$ be a monotone pair of nonnegative vectors. Then, we have
$$
\sum a\!\cdot\! X(b) \le {\rm rc}(X)\sum X(a\!\cdot\! b)
$$
and ${\rm rc}(X)$ is the best possible constant not depending  on $(a,b)$.
}

\vskip 10pt\noindent
{\bf Lemma 1.4.} {\it Let $X$ be a real  sum-symmetric matrix and let $a$ and  $b$ be  vectors. 
If $(a,b)$ is monotone,
$$
\sum a\!\cdot\! X(b) \le \sum X(a\!\cdot\! b).
$$
If $(a,b)$ is antimonotone, the reverse inequality holds.
}

\vskip 10pt\noindent
{\bf Proof of Lemma 1.3.} Since ${\rm rc} (X)={\rm rc} (PXP^T)$ for any permutation matrix $P$, and since
$$
\sum a\!\cdot\! X(b)=\sum Pa\!\cdot\! PXP^T(Pb) \quad{\rm and}\quad \sum X(a\!\cdot\! b)=\sum PXP^T(P(a)\!\cdot\!P(b)),
$$
we may assume that both the components of $a$ and $b$ are arranged in increasing order.
Let $e_1$ be the vector with all components 1, $e_2$ the vector with the first component 0 and all the others 1, ....\,, and  $e_n$ the vector with last component 1 and all the others 0. There exist nonnegative scalars $\alpha_j$ and $\beta_j$ such that
$$
a=\sum_{1\le j\le n} \alpha_j e_j \quad{\rm and}\quad b=\sum_{1\le j\le n} \beta_j e_j.
$$
By linearity, it suffices to show
\begin{equation}
\sum e_j\!\cdot\! X(e_k) \le {\rm rc}(X)\sum X(e_j\!\cdot\! e_k)
\end{equation}
for all indexes $j$, $k$. We distinguish two cases.

If $j\le k$, then $X(e_j \!\cdot\! e_k)=X(e_k)$ and the inequality is obvious since $\sum e_j\!\cdot\! X(e_k) \le \sum X(e_k)$.
If $j> k$, then using the definition of ${\rm rc}(X)$,
\begin{align*}
\sum e_j\!\cdot\! X(e_k) &=\sum_{l=j}^n \sum_{m=k}^n x_{l,m}\cr
&\le\sum_{l=j}^n \sum_{m=1}^n x_{l,m}\cr
&\le\sum_{l=j}^n {\rm rc}(X) \sum_{m=1}^n x_{m,l}\cr
&={\rm rc}(X) \sum_{m=1}^n \sum_{l=j}^n x_{m,l}\cr
&= {\rm rc}(X)\sum X(e_j) \cr
&={\rm rc}(X)\sum X(e_j\!\cdot\! e_k).
\end{align*}
Hence (11) holds and the main part of Lemma 1.3 is proved. To see that this inequality is sharp, consider   a vector $u$ of the canonical basis correponding with an index $i_0$ such that 
$$
{\rm rc}(X)=\frac{\sum_k x_{i_0,k}}{ \sum_k x_{k,i_0}}.
$$
Then
$$
{\rm rc}(X)=\frac{\sum X^T(u)}{\sum X(u)}=\frac {\sum uX(e_1)}{\sum X(ue_1)}
$$
and  $(u,e_1)$ is monotone. This completes the proof. \qquad $\Box$

\vskip 10pt\noindent
{\bf Proof of Lemma 1.4.} Since $(a,b)$ is monotone iff $(a,-b)$ is antimonotone, it suffices to consider the first case.
Fix a monotone pair $(a,b)$ and  a constant $\gamma$. Then $(a,b)$ satisfies to  the lemma iff the same holds for $(a+\gamma e_1,b+\gamma e_1)$. Hence we may suppose that $(a,b)$ is nonnegative and we apply Lemma 1.3 with ${\rm rc} (X)=1$. \qquad $\Box$

\vskip 10pt
Let us show how Theorem 1.1 (and similarly Theorem 1.2) follows from Lemma 1.4. Since $(A,B)$ is monotone, we may assume that $A$ and $B$ are diagonal, 
$A={\rm diag}(\alpha_1,\dots,\alpha_n)$ and $B={\rm diag}(\beta_1,\dots,\beta_n)$. Let $X$ be the  sum-symmetric matrix $X=(x_{i,j})$ with $x_{i,j}=|z_{i,j}|^2$ and oberve that
$$
\Vert AZB\Vert_2^2=\sum a\!\cdot\! X(b)\quad {\rm and}\quad \Vert ZAB\Vert_2^2=\sum X(a\!\cdot\! b)
$$
where $a=(\alpha_1^2\dots\alpha_n^2)^T$ and $b=(\beta_1^2\dots\beta_n^2)^T$ form a monotone or an antimonotone pair of vectors.

\vskip 10pt
The following result [4] is another extension of (10). We omit the proof since it is contained in Theorem 2.1 below.

\vskip 10pt\noindent
{\bf Theorem 1.5.} {\it Let $A,\, B\ge 0$ with $(A,B)$ monotone and let $E$ be a projection. Then, there exists a unitary $U$ such that
\begin{equation*}
\vert AEB\vert \le U \vert ABE\vert U^*.
\end{equation*} 
}

From this result we derived several eigenvalues inequalities and a determinantal Chebyshev type inequality involving compressions: {\it Let $(A,B)$ be monotone positive and let ${\cal E}$ be a subspace. Then,}
$$
\det A_{\cal E}\cdot\det B_{\cal E}\le \det (AB)_{\cal E}.
$$
Here $A_{\cal E}$ denotes the compression of $A$ onto ${\cal E}$. When ${\cal E}$ has codimension 1, we showed the reverse inequality for antimonotone positive pairs $(A,B)$,
$$
\det A_{\cal E}\cdot\det B_{\cal E}\ge \det (AB)_{\cal E}.
$$
Of course, this also holds for 1-dimensional subspaces as a restatement of (10). Hence we raised the following question:

\vskip 10pt\noindent
{\bf Problem 1.6.} Does the above determinantal inequality for antimonotone pairs hold for all subspaces?

\vskip 10pt\noindent

We mention another open problem. We showed  [3] that Theorem 1.2 can not be extended to Schatten $p$-norms when $p>2$ by giving counterexamples in dimension 3. But the following is still open:

\vskip 10pt\noindent
{\bf Problem 1.7.} Does Theorem 1.1 hold for Schatten $p$-norms, $1\le p<2$? In particular for the Trace norm?

\vskip 15 pt\noindent
{\large\it 1.1 Gruss type inequalities for the trace}
\vskip 10pt
In connection with Theorem 1.2 we have
the following two results which are  Gruss type inequalities for the trace. Letting $Z$ be a rank one projection in the first result and assuming $AB=BA$ we get the classical Gruss inequalities (3), (4).

\vskip 10pt\noindent
{\bf Proposition 1.8.} {\it
For $Z\ge0$, Hermitian $A$ with extremal eigenvalues $p$ and $q$ $(p\ge q)$ and Hermitian $B$ with extremal eigenvalues $r$ and $s$ $(r\ge s)$, 
$$
\left|  {\rm Tr\,}Z^2AB-{\rm Tr\,}ZAZB \right| \le \frac{1}{4}(p-q)(r-s){\rm Tr\,}Z^2.
$$
In particular,
$$
{\rm Tr\,}Z^2A^2-{\rm Tr\,}(ZA)^2 \le \frac{(p-q)^2}{4}{\rm Tr\,}Z^2.
$$ 
}

\vskip 10pt\noindent
{\bf Proof.} Note that 
$$
{\rm Tr\,}Z^2AB = \overline{{\rm Tr}\,(Z^2AB)^*}= \overline{{\rm Tr}\,Z^2BA}
$$
and similarly
$$
{\rm Tr\,}ZAZB = \overline{{\rm Tr}\,(ZBZA)}.
$$
 Since $(A,A)$ is monotone, Theorem 1.2  shows that the map
$$
(A,B)\longrightarrow  {\rm Tr\,}Z^2AB-{\rm Tr\,}ZAZB
$$
is a complex valued semi-inner product on the real vector space of Hermitian operators. The Cauchy-Schwarz inequality for this semi-inner product then shows that it suffices to prove the second inequality of our theorem. Let $Z=\sum_i z_i e_i\otimes e_i$ be the canonical expansion of $Z$. Since the Frobenius norm of a matrix is less than the $l^2$-norm of its diagonal, we have 
\begin{align*}
{\rm Tr\,}Z^2A^2-{\rm Tr\,}(ZA)^2 &= {\rm Tr\,}Z^2A^2- \Vert Z^{1/2}AZ^{1/2} \Vert_2^2 \\
&\le \sum_i z_i^2\langle e_i, A^2e_i\rangle - \sum_i (z_i\langle e_i, Ae_i\rangle)^2 \\
&\le \frac{(p-q)^2}{4}\sum_i z_i^2 = \frac{(p-q)^2}{4}{\rm Tr\,}Z^2
\end{align*}
by using the classical inequality (4). \qquad $\Box$

\vskip 10pt
Letting $Z$ be a rank one projection we recapture an inequality pointed out by M.\ Fujii et al. [15],
$$
\left| \langle h, ABh\rangle - \langle h, Ah\rangle\langle h, Bh\rangle \right|\le \frac{1}{4}(p-q)(r-s)
$$
for all norm one vectors $h$. They called it the Variance-covariance Inequality. By using the GNS construction, this can be formulated in the $C^*$-algebra framework: Given positive elements $a$, $b$ with spectra in $[p,q]$ and $[r,s]$ respectively,
$$
\left| \varphi(ab) - \varphi(a)\varphi(b) \right|\le \frac{1}{4}(p-q)(r-s)
$$
for all states $\varphi$.

\vskip 10pt\noindent
{\bf Proposition 1.9.} {\it
For normal $Z$, Hermitian $A$ with extremal eigenvalues $p$ and $q$ $(p\ge q)$ and Hermitian $B$ with extremal eigenvalues $r$ and $s$ $(r\ge s)$,
$$
\left|  {\rm Tr\,}|Z|^2AB-{\rm Tr\,}Z^*AZB \right| \le \frac{1}{2}(p-q)(r-s){\rm Tr\,}|Z|^2.
$$
}

 \noindent 
{\bf Proof.} Theorem 1.2 shows that the map
$$
(A,B)\longrightarrow  {\rm Tr\,}|Z|^2AB-{\rm Tr\,}Z^*AZB
$$
is a semi-inner product on the space of Hermitian operators. Hence it suffices to consider the case $A=B$:

Let 
$$
\tilde{Z}=
\begin{pmatrix}
 0&Z^* \\
 Z&0
\end{pmatrix}
 \qquad {\rm and} \qquad
 \tilde{A}=
\begin{pmatrix}
 A&0 \\
 0&A
\end{pmatrix}
$$
and observe that
$$
{\rm Tr\,}|Z|^2A^2-{\rm Tr\,}Z^*AZA =\frac{1}{2}\{\,{\rm Tr\,}\tilde{Z}^2\tilde{A}^2-{\rm Tr\,}(\tilde{Z}\tilde{A})^2\,\}
$$
and
$$
{\rm Tr\,}|Z|^2=\frac{1}{2}{\rm Tr\,}\tilde{Z}^2.
$$
Consequently, $\tilde{Z}$ being Hermitian, we may assume so is $Z$. Replacing if necessary $A$ by $A+qI$, we may also assume $A\ge 0$. We compute in respect with an orthonormal basis of eigenvectors for $A=\sum_i a_i e_i\otimes e_i$,
\begin{align*}
{\rm Tr\,}Z^2A^2-{\rm Tr\,}(ZA)^2 &={\rm Tr\,}Z^2A^2 -\Vert A^{1/2}ZA^{1/2}\Vert_2^2 \\
&=\sum_{i,j}a_i^2|z_{i,j}|^2-\sum_{i,j}a_ia_j|z_{i,j}|^2 \\
&=\sum_{i<j}(a_i-a_j)^2|z_{i,j}|^2 \\
&\le \frac{(p-q)^2}{2}{\rm Tr\,}Z^2
\end{align*}
and the proof is complete. \qquad $\Box$

\vskip 20 pt\noindent
{\large\bf 2. Matrix Kantorovich  inequalities}
\vskip 10pt
In view of Theorem 1.5 involving projections, we tried to obtain an extension for all positive operators. We obtained a hybrid Chebyshev/Kantorovich result:

\vskip 10pt\noindent
{\bf Theorem 2.1.} {\it Let  $A,\,B\ge 0$ with $(A,B)$ monotone and let $Z\ge0$ with its largest and smallest nonzero eigenvalues  $a$ and $b$. Then, there exists a unitary  $U$ such that
\begin{equation*}
|AZB|\le \,\frac{a+b}{2\sqrt{ab}}\, U|ZAB|U^*.
\end{equation*}
}
\noindent
Since for a projection $Z$, we have $a=b=1$, Theorem 2.1 contains Theorem 1.5. 

\vskip 10pt
 We recall that the  inequality of Theorem 2.1 is equivalent to:
$$
\mu_j(AZB) \le \mu_j(ZAB)
$$
for all $j=1,\dots$, where $\{\mu_j(\cdot)\}$ stand for the singular values arranged in decreasing order with their multiplicities [2, p.\ 74].

\vskip 10pt
We need some lemmas. First we state the Kantorovich inequality (7) again and we give a matrix proof.

\vskip 10pt\noindent
{\bf Lemma 2.2.} {\it Let  $Z>0$ with extremal eigenvalues  $a$ and $b$. Then, for every norm one vector $h$,
$$
\Vert Zh\Vert \le \,\frac{a+b}{2\sqrt{ab}}\, \langle h,Zh\rangle.
$$
}
 
\vskip 10pt\noindent
{\bf Proof.} Let ${\cal E}$ be any subspace of ${\cal H}$ and let $a'$ and $b'$ be the extremal eigenvalues of $Z_{\cal E}$. Then $a\ge a'\ge b'\ge b$ and, setting $t=\sqrt{a/b}$, $t'=\sqrt{a'/b'}$, we have $t\ge t'\ge1$. Since $t\longrightarrow t+1/t$ increases on $[1,\infty)$ and
$$
\frac{a+b}{2\sqrt{ab}}\,=\frac{1}{2}\left( t+ \frac{1}{t}\right), \qquad
 \frac{a'+b'}{2\sqrt{a'b'}}\, =\frac{1}{2}\left( t'+ \frac{1}{t'}\right),
$$
we infer
$$
\frac{a+b}{2\sqrt{ab}}\,\ge \frac{a'+b'}{2\sqrt{a'b'}}\,.
$$
Therefore, it suffices to prove the lemma for $Z_{\cal E}$ with ${\cal E}={\rm span}\{h,Zh\}$. Hence, we may assume $\dim{\cal H}=2$, $Z=ae_1\otimes e_1+be_2\otimes e_2$ and $h=xe_1 +(\sqrt{1-x^2})e_2$. Setting
$x^2=y$ we have
$$
\frac{||Zh||}{\langle h, Zh\rangle}=\frac{\sqrt{a^2y+b^2(1-y)}}{ay+b(1-y)}.
$$
The right hand side attains its maximum on $[0,1]$ at $y=b/(a+b)$,
and then
$$
\frac{||Zh||}{\langle h, Zh\rangle}=\,\frac{a+b}{2\sqrt{ab}}
$$
proving the lemma.\qquad $\Box$

\vskip 10pt\noindent
Lemma 2.2 can be extended as an inequality involving the operator norm $\Vert\cdot\Vert_{\infty}$ and the spectral radius $\rho(\cdot)$. Indeed, letting $A=h\otimes h$ in the next lemma, we get Lemma 2.2. 

\vskip 10pt\noindent
{\bf Lemma 2.3.} {\it
 For $A\ge0$ and  $Z>0$ with extremal eigenvalues  $a$ and $b$,
\begin{equation*}
\Vert AZ\Vert_{\infty}\le \,\frac{a+b}{2\sqrt{ab}}\, \rho(AZ).
\end{equation*}
}

\vskip 5pt\noindent
{\bf Proof.} There exists a rank one projection $F$ such that, letting $f$ be a unit vector in the range of $A^{1/2}F$,
\begin{align*}
\Vert AZ\Vert_{\infty}=\Vert ZA\Vert_{\infty} 
=\Vert ZAF\Vert_{\infty}&=\Vert ZA^{1/2}(f\otimes f)A^{1/2}F\Vert_{\infty} \\
&\le \Vert ZA^{1/2}(f\otimes f)A^{1/2}\Vert_{\infty} 
= \Vert A^{1/2} f\Vert^2 \Vert Z \frac{A^{1/2}f}{\Vert A^{1/2} f\Vert} \Vert.
\end{align*}
Hence
\begin{align*}
\Vert AZ\Vert_{\infty} \le \,\frac{a+b}{2\sqrt{ab}}\, \langle f, A^{1/2}ZA^{1/2}f\rangle
\le \,\frac{a+b}{2\sqrt{ab}}\, \rho(A^{1/2}ZA^{1/2})= \,\frac{a+b}{2\sqrt{ab}}\,\rho(AZ)
\end{align*}
by using Lemma 2.2 with $h=A^{1/2}f/\Vert A^{1/2} f\Vert $. \qquad $\Box$

\vskip 10pt\noindent
From Lemma 2.3 one may derive a sharp operator inequality:
 \vskip 10pt\noindent
{\bf Lemma 2.4.} {\it Let $0\le A\le I$ and let $Z>0$ with extremal eigenvalues  $a$ and $b$. Then, 
$$
AZA \le \,\frac{(a+b)^2}{4ab}\, Z.
$$
}

\noindent
{\bf Proof.} The claim is equivalent to the operator norm inequalities
$$
\Vert Z^{-1/2}AZAZ^{-1/2}\Vert_{\infty} \le \,\frac{(a+b)^2}{4ab}
$$
or
$$
\Vert Z^{-1/2}AZ^{1/2}\Vert_{\infty} \le \,\frac{a+b}{2\sqrt{ab}}.
$$
But the previous lemma entails
\begin{align*}
\Vert Z^{-1/2}AZ^{1/2}\Vert_{\infty}&=\Vert Z^{-1/2}AZ^{-1/2}Z\Vert_{\infty} \\
&\le \,\frac{a+b}{2\sqrt{ab}}\,\rho(Z^{-1/2}AZ^{-1/2}Z) \\
& =\,\frac{a+b}{2\sqrt{ab}}\,\Vert A\Vert_{\infty} \\
& \le\,\frac{a+b}{2\sqrt{ab}},
\end{align*}
hence, the result holds. \qquad $\Box$

\vskip 10pt\noindent
{\bf Proof of Theorem 2.1.} We will use Lemma 2.4 and the following operator norm inequality
\begin{equation}
\Vert AEB\Vert_{\infty} \le \Vert ABE\Vert_{\infty}
\end{equation}
for all projections $E$. This inequality was derived [3] from Theorem 1.1  and is the starting point and a special case of Theorem 1.5. In fact (12) is a consequence of (10).
Indeed there exist unit vectors $f$ and $h$ with $h=Eh$ such that
$$
\Vert AEB\Vert_{\infty}=\Vert AEBf\Vert=\Vert A(h\otimes h)Bf\Vert\le \Vert A(h\otimes h)B\Vert_{\infty}=\Vert Ah\Vert \Vert Bh\Vert
$$
so that using (10)
$$
\Vert AEB\Vert_{\infty}\le \Vert ABh\Vert \le\Vert ABE\Vert_{\infty}.
$$

We denote by ${\rm supp}(X)$ the support projection of an operator $X$, i.e., the smallest projection $S$
such that $X=XS$.

By the minimax principle, for every projection $F$, ${\rm corank} F=k-1$,
\begin{align}
\mu_k(AZB)&\le \Vert AZBF\Vert_{\infty}  \\
&= \Vert AZ^{1/2}EZ^{1/2}BF\Vert_{\infty} \notag \\
&\le \Vert AZ^{1/2}EZ^{1/2}B\Vert_{\infty} \notag
\end{align}
where $E$ is the projection onto the range of $Z^{1/2}BF$. Note that there exists a rank one projection $P$, $P\le E$,  such that
$$
\mu_k(AZB)\le \Vert AZ^{1/2}PZ^{1/2}B\Vert_{\infty}.
$$
Indeed, let $h$ be a norm one vector such that
$$
\Vert AZ^{1/2}EZ^{1/2}B\Vert_{\infty}=\Vert AZ^{1/2}EZ^{1/2}Bh\Vert
$$
and let $P$ be the projection onto ${\rm span}\{EZ^{1/2}Bh\}$. Since $Z^{1/2}PZ^{1/2}$
has rank one, and hence is a scalar multiple of a projection, (12) entails
\begin{equation*}
\mu_k(AZB)\le \Vert Z^{1/2}PZ^{1/2}AB\Vert_{\infty}. 
\end{equation*}
We may  choose $F$ in (13) in order to obtain any projection $G\ge {\rm supp}(EZ^{1/2}AB)$, ${\rm corank} G=k-1$. Since
$$
{\rm supp}(PZ^{1/2}AB)\le {\rm supp}(EZ^{1/2}AB)\le G,
$$
we infer
\begin{equation*}
\mu_k(AZB)\le \Vert Z^{1/2}PZ^{1/2}ABG\Vert_{\infty}.
\end{equation*}
Consequently, using Lemma 2.4 with $Z$ and $PZP$,
\begin{align*}
\mu_k(AZB)&= \Vert GABZ^{1/2}PZPZ^{1/2}ABG\Vert_{\infty}^{1/2} \\
&\le \,\frac{a+b}{2\sqrt{ab}}\,\Vert ZABG\Vert_{\infty}. 
\end{align*}
Since we may choose $G$ so that $\Vert ZABG\Vert_{\infty}=\mu_k(ZAB)$, the proof is complete.
 \qquad $\Box$
 
 \vskip 10pt
 Under an additional invertibility assumption on $Z$, Theorem 2.1 can be reversed: 
 \vskip 10pt\noindent
{\bf Theorem 2.5.} {\it Let  $A,\,B\ge 0$ with $(A,B)$ monotone and let $Z>0$ with extremal eigenvalues  $a$ and $b$. Then, there exists a unitary  $V$ such that
\begin{equation*}
|ZAB|\le \,\frac{a+b}{2\sqrt{ab}}\, V|AZB|V^*
\end{equation*}
}

\vskip 10pt\noindent
{\bf Proof.} By a limit argument we may assume that both  $A$ and $B$ are invertible. Hence, taking inverses in Theorem 2.1 considered as singular values inequalities, we obtain a unitary $W$ (actually we can take $W=U$ since $t\longrightarrow t^{-1}$ is operator decreasing) such that
$$
|AZB|^{-1}\ge \frac{2\sqrt{ab}}{a+b}\,  W|ZAB|^{-1}W^* 
$$
hence, using $|X|^{-1}=(X^*X)^{-1/2}=(X^{-1}X^{*-1})^{1/2}=|X^{*-1}|$ for all invertibles $X$,
$$
|B^{-1}Z^{-1}A^{-1}|\ge \frac{2\sqrt{ab}}{a+b}\,  W|A^{-1}B^{-1}Z^{-1}|W^*.
$$
Then observe that we can replace $Z^{-1}$ by $Z$ since 
$$
\,\frac{a+b}{2\sqrt{ab}}\,=\,\frac{a^{-1}+b^{-1}}{2\sqrt{a^{-1}b^{-1}}}.
$$
As the correspondence between an invertible monotone pair and its inverse is onto, Theorem 2.5 holds. \quad $\Box$ 

\vskip 10pt
Let $X$ with real eigenvalues and  denote by $\lambda_k(X)$, $k=1,\,2,\dots$, the eigenvalues of $X$ arranged in decreasing order with their multiplicities. Replacing $A$ and $B$ by $A^{1/2}$ in Theorems 2.1, 2.5 we get:

\vskip 10pt\noindent
{\bf Corollay 2.6.} {\it Let  $A\ge 0$ and let $Z>0$ with extremal eigenvalues  $a$ and $b$. Then,
for all $k$,
\begin{equation*}
\frac{2\sqrt{ab}}{a+b}\, \lambda_k (AZ)\le \mu_k(AZ)\le \,\frac{a+b}{2\sqrt{ab}}\, \lambda_k(AZ).
\end{equation*}
}
\vskip 10pt\noindent
Note that Corollary 2.6 contains Lemma 2.3, hence Lemma 2.2. 

\vskip 10pt
By replacing in Theorem 2.5 $A$ and $B$ by a rank one projection $h\otimes h$ we recapture the Kantorovich inequality of
Lemma 2.2. This shows that Theorem 2.5 is sharp and, since they are equivalent, also Theorem 2.1 (see [8] for more details). Similarly to Theorem 2.5, the next theorem is also a sharp inequality extending Lemma 2.2.

\vskip 10pt\noindent
{\bf Theorem 2.7.} {\it Let  $A$, $B$ such that $AB\ge0$ and let $Z>0$ with extremal eigenvalues  $a$ and $b$. Then, for all symmetric norms,
$$
\Vert ZAB\Vert \le \,\frac{a+b}{2\sqrt{ab}}\, \Vert BZA\Vert.
$$
}

\vskip 10pt\noindent
As in the special case of the operator norm (1), a basic rearrangement inequality for general symmetric norms claims that
\begin{equation}
\Vert AB\Vert\le\Vert BA\Vert
\end{equation}
whenever the product $AB$ is normal. Thus, when $AB\ge 0$ Theorem 2.7 is a  generalization of (14).
Let us give a proof of (14).
First for all symmetric norms and all partitionned matrices,
$$
\left|\left| \begin{pmatrix}
A&0 \\ 0&B
\end{pmatrix} \right|\right|
\le
\left|\left| \begin{pmatrix}
A&R \\ S&B
\end{pmatrix} \right|\right|,
$$
indeed, the left hand side is the mean of two unitary conguences of the right hand side,
$$
 \begin{pmatrix}
A&0 \\ 0&B
\end{pmatrix} 
 = \frac{1}{2}
 \begin{pmatrix}
A&R \\ S&B
\end{pmatrix}
+ 
\frac{1}{2}
 \begin{pmatrix}
I&0 \\ 0&-I
\end{pmatrix}
 \begin{pmatrix}
A&R \\ S&B
\end{pmatrix} 
\begin{pmatrix}
I&0 \\ 0&-I
\end{pmatrix}.
$$
By repetiton of this argument we see that symmetric norms of any  matrix  $X$ are greater than those of its diagonal,
\begin{equation}
\Vert {\rm diag}(X)\Vert \le \Vert X\Vert.
\end{equation}
This inequality is quite important. Applying (15) to $X=BA$ with $AB$ normal we deduce, by writing $X$ in a triangular form, that
$$
\Vert AB\Vert =\Vert {\rm diag}(BA)\Vert \le \Vert BA\Vert.
$$
Therefore (14) holds.

\vskip 10pt\noindent
{\bf Proof of Theorem 2.7.} Using Corollary 2.6 we have
\begin{align*}
\Vert ZAB\Vert = \Vert {\rm diag}(\mu_k(ZAB))\Vert &\le \,\frac{a+b}{2\sqrt{ab}}\,\Vert {\rm diag}(\lambda_k(ZAB))\Vert \\
&=\,\frac{a+b}{2\sqrt{ab}}\,\Vert {\rm diag}(\lambda_k(BZA))\Vert\le \,\frac{a+b}{2\sqrt{ab}}\,\Vert BZA\Vert
\end{align*}
where the last inequality follows from (15) applied to $BZA$ in a triangular form. \qquad $\Box$

\vskip 10pt\noindent
{\bf Remarks.}  Starting from Lemma 2.2, we first proved  Theorem 2.7  in [5] by using Ky Fan dominance principle (see the next section). As applications we then derived the above Lemmas 2.3 and 2.4. Theorem 2.1 had been proved later [8]. In some sens, the presentation given here, which starts from the earlier Theorem 1.5, is more natural. From Lemma 2.4 we also derived:

\noindent
 (Mond-Pe$\check{{\rm c}}$ari\'c [22]) {\it Let $Z>0$ with extremal eigenvalues  $a$ and $b$. Then, for every subspace ${\cal E}$, 
$$
(Z_{\cal E})^{-1} \ge \,\frac{4ab}{(a+b)^2}\,(Z^{-1})_{\cal E}.
$$
}
A similar compression inequality holds for others operator convex functions [5]. Mond-Pe$\check{{\rm c}}$ari\'c's result is clearly an extension of the original Kantorovich inequality (6), (7).

\vskip 10pt\noindent
{\bf Corollary 2.8.} {\it Let  $A,\,B>0$ with $AB=BA$ and $pI\ge AB^{-1}\ge qI$ for some  $p,q>0$. Then, for all $Z\ge 0$ and all symmetric norms
$$
\Vert AZB\Vert \le \,\frac{p+q}{2\sqrt{pq}}\, \Vert ZAB\Vert.
$$
}

\vskip 5pt\noindent
{\bf Proof.} Write $AZB=AZA(A^{-1}B)$ and apply Theorem 2.7 with $A^{-1}B$ instead of $Z$. \qquad $\Box$

 \vskip 15 pt\noindent
{\large\it 2.1. Rearrangement inequalities for sequences}
\vskip 10pt
Corollary 2.8 can not be extended to normal operators $Z$, except in the case of the trace norm. This observation leaded to establish [9] the following reverse inequality to the most basic rearrangement inequality (9). Recall that down arrows mean nonincreasing rearrangements.

\vskip 10pt\noindent
{\bf Theorem 2.9.} {\it Let $\{a_i\}_{i=1}^n$ and $\{b_i\}_{i=1}^n$ be $n$-tuples of positive numbers with
$$
 p\ge \frac{a_i}{b_i}\ge q,\quad i=1,\dots,n,
$$
for some $p,q>0$. Then,
$$
\sum_{i=1}^n a_i^{\downarrow}b_i^{\downarrow} \le \frac{p+q}{2\sqrt{pq}}\, \sum_{i=1}^n a_ib_i.
$$
}

\noindent
{\bf Proof.} Introduce the diagonal matrices $A={\rm diag}(a_i)$ and $B={\rm diag}(b_i)$ and observe that, $\Vert\cdot\Vert_1$ standing for the trace norm,
$$
\sum_{i=1}^n a_ib_i=\Vert AB\Vert_1
$$
and 
$$
\sum_{i=1}^n a_i^{\downarrow}b_i^{\downarrow}= \Vert AVB\Vert_1
$$
for some permutation matrix $V$. Hence we have to show that
$$
\Vert AVB\Vert_1 \le \frac{p+q}{2\sqrt{pq}}\, \Vert AB\Vert_1
$$
To this end consider the spectral representation $V=\sum_i v_i h_i\otimes h_i$ where $v_i$ are the eigenvalues and $h_i$
the corresponding unit eigenvectors. We have 
\begin{align*}
\Vert AVB\Vert_1 &\le \sum_{i=1}^n \Vert A\cdot  v_i h_i\otimes h_i\cdot B \Vert_1 \\
&= \sum_{i=1}^n \Vert A h_i\Vert \, \Vert Bh_i \Vert \\
&\le \frac{p+q}{2\sqrt{pq}} \,\sum_{i=1}^n\langle Ah_i, Bh_i\rangle \\
&= \frac{p+q}{2\sqrt{pq}}\, \sum_{i=1}^n\langle h_i, ABh_i\rangle \\
&= \frac{p+q}{2\sqrt{pq}} \,\Vert AB\Vert_1
\end{align*}
where we have used the triangle inequality for the trace norm and  Lemma 2.10 below. \qquad $\Box$

\vskip 10pt\noindent
{\bf Lemma 2.10} {\it Let  $A,\,B>0$ with $AB=BA$ and $pI\ge AB^{-1}\ge qI$ for some  $p,q>0$. Then, for every  vector $h$,
$$
\Vert Ah\Vert \,\Vert Bh\Vert\le \,\frac{p+q}{2\sqrt{pq}}\, \langle Ah,Bh\rangle.
$$
}

\vskip 10pt\noindent
{\bf Proof.} Write $h=B^{-1}f$ and apply Lemma 2.2; or apply Corollary 2.8 with $Z=h\otimes h$. \qquad $\Box$

\vskip 10pt\noindent
{\bf Remark.} Lemma 2.10 extends Lemma 2.2 and is nothing less but  of Cassel's Inequality:
\vskip 10pt\noindent
{\it Cassel's inequality}. For nonnegative $n$-tuples $\{a_i\}_{i=1}^n$, $\{b_i\}_{i=1}^n$ and $\{w_i\}_{i=1}^n$ with
$$
 p\ge \frac{a_i}{b_i}\ge q,\quad i=1,\dots,n,
$$
for some $p,\,q>0$\,; it holds that
$$
\left(\sum_{i=1}^n w_ia_i^2\right)^{1/2}\left(\sum_{i=1}^n w_ib_i^2\right)^{1/2} \le \frac{p+q}{2\sqrt{pq}}\, \sum_{i=1}^n w_ia_ib_i.
$$
Of course it is a reverse inequality to the Cauchy-Schwarz inequality. To obtain it from Lemma 2.9, one just
takes $A={\rm diag}(a_1,\dots,a_n)$, $B={\rm diag}(b_1,\dots,b_n)$ and $h=(\sqrt{w_1},\dots,\sqrt{w_n})$. If one let $a=(a_1,\dots,a_n)$ and $b=(b_1,\dots,b_n)$ then Cassel's inequality can be written 
\begin{equation*}
\Vert a\Vert\, \Vert b\Vert \le \frac{p+q}{2\sqrt{pq}}\,\langle a,b\rangle 
\end{equation*}
for a suitable inner product $\langle \cdot,\cdot \rangle $.
It is then natural to search for conditions on $a$, $b$ ensuring that the above inequality remains valid with $Ua$, $Ub$
for all orthogonal matrices $U$. This motivates a remarkable extension of Cassel's inequality:

\vskip 10pt\noindent
{\it Dragomir's inequality}. For real vectors $a$, $b$ such that $\langle a-qb, pb-a\rangle\ge 0$ for some scalars $p,\,q$ with $pq>0$, inequality (1) holds.

\vskip 10pt\noindent
 Dragomir's inequality admits a version for complex vectors. For these inequalities see [10], [11], [12].

\vskip 10pt\noindent
{\bf Remark.} In  [9] we also investigate reverse additive inequalities to (9). This setting is less clear.
In general reverse additive type inequalities are more difficult than multiplicative ones. The story of Ozeki's inequality, a reverse additive inequality to Cauchy-Schwarz's inequality illustrates that [19].

\vskip 20pt\noindent
{\large\bf 3. Generalized Kantorovich inequalities}
\vskip 10pt

In [1] (see also [2, pp.\ 258, 285]) Araki showed a trace inequality which entails the following inequality for symmetric norms:

\vskip 10pt\noindent
{\bf Theorem 3.1.} {\it Let  $A\ge 0$, $Z\ge 0$ and $p>1$.  Then, for every symmetric norm, 
\begin{equation*}
\Vert (AZA)^p\Vert \le  \Vert A^pZ^pA^p\Vert.
\end{equation*}
For $0<p<1$, the above inequality is reversed.
}

\vskip 10pt
If we take  a rank one projection $A=h\otimes h$, $\Vert h\Vert=1$, then Araki's inequality  reduces to Jensen's inequality for $t\longrightarrow t^p$,
\begin{equation}
\langle h, Zh\rangle^p \le \langle h,Z^ph\rangle.
\end{equation}
This inequality admits a reverse inequality. Ky Fan [20]  introduced the following constant, for $a$,\,$b>0$ and integers $p$,
$$
K(a,b,p)=\frac{a^pb-ab^p}{(p-1)(a-b)}\left(\frac{p-1}{p}\frac{a^p-b^p}{a^pb-ab^p}\right)^p.
$$
Furuta extended it to all real numbers (see for instance [17], [18])
and showed the sharp reverse inequality of (16): {\it If $Z>0$ have extremal eigenvalues $a$ and $b$, then
\begin{equation}
\langle h, Z^ph\rangle \le K(a,b,p)\langle h,Zh\rangle^p
\end{equation}
for $p>1$ and $p<0$.}

\vskip 10pt\noindent
In a recent paper [14], Fujii-Seo-Tominaga extended (17) to an operator norm inequality: {\it For $A\ge 0$,   $Z>0$ with extremal eigenvalues  $a$ and $b$, and $p>1$,
$$
\Vert A^pZ^pA^p\Vert \le  K(a,b,p)\Vert (AZA)^p\Vert_{\infty}.
$$
}
Inspired by this result, we showed in [7]:

\vskip 10pt\noindent
{\bf Theorem 3.2.} {\it Let  $A\ge 0$ and let $Z>0$ with extremal eigenvalues  $a$ and $b$. Then, for every $p>1$, there exist unitaries $U$,  $V$ such that 
\begin{equation*}
\frac{1}{ K(a,b,p)}\,U(AZA)^pU^*\le A^pZ^pA^p \le  K(a,b,p)V(AZA)^pV^*.
\end{equation*}
The Ky Fan constant $K(a,b,p)$ and its inverse are optimal.
}

\vskip 10pt\noindent
For $p=2$, Theorem 3.2 is a reformulation of Corollary 2.6.

 Furuta introduced another constant depending on reals $a$, $b$ and $p>1$
$$
C(a,b,p)=(p-1)\left( \frac{a^p-b^p}{p(a-b)}\right)^{p/(p-1)} + \frac{ab^p-ba^p}{a-b}
$$
in order to obtain
\begin{equation}
\langle h, Z^ph\rangle - \langle h, Zh\rangle^p \le C(a,b,p)
\end{equation}
for unit vectors $h$ and $Z\ge0$ with extremal eigenvalues $a$ and $b$ (see [24, Theorem C]). Equivalently,
$$
C(a,b,p)=\max \{ \int_{\Omega} f^p\,{\rm d}\mu -\left(\int_{\Omega} f\,{\rm d}\right)^p \}
$$
where the maximum runs over all measurable functions $f$, $a\ge f(t)\ge b$, on probabilized space $(\Omega,\mu)$. Hence
$a\ge a'\ge b' \ge b \Rightarrow C(a,b,p)\ge C(a',b',p)$.

Of course (18) generalizes the quadratic case (4) and $C(a,b,2)=(a-b)^2/4$. Simplified proofs are given in [14] by using the Mond-Pe$\check{{\rm c}}$ari\'c method. It is also possible to prove it by reduction to the $2\times 2$ matrix case, in a similar way of Lemma 2.2.

Furuta's constant allows us to extend the second inequality of Proposition 1.8 (in which $p=2$):

\vskip 10pt\noindent
{\bf Lemma 3.3.} {\it Let $A\ge0$ and let  $Z\ge0$ with extremal eigenvalues $a$ and $b$. Then, for all $p>1$,
$$
{\rm Tr}\, A^pZ^pA^p -  {\rm Tr}\,(AZA)^p \le C(a,b,p)  {\rm Tr}\,A^{2p}.
$$
}
This trace inequality can be extended to all symmetric norms:

\vskip 10pt\noindent
{\bf Theorem 3.4.} {\it Let $A\ge0$ and let  $Z\ge0$ with extremal eigenvalues $a$ and $b$. Then for all symmetric norms and all $p>1$,
$$
\Vert A^pZ^pA^p\Vert - \Vert (AZA)^p\Vert \le C(a,b,p) \Vert A^{2p}\Vert.
$$
}

\noindent
Note that letting $A$ be a rank one projection either in the lemma or the theorem, we recapture  inequality (18).

\vskip 10pt
 We will use the Ky Fan dominance principle: {\it  $\Vert A\Vert \le \Vert B\Vert$ for all symmetric norms iff 
$\Vert A\Vert_{(k)} \le \Vert B\Vert_{(k)}$ for all Ky Fan $k$-norms}. By definition $\Vert A\Vert_{(k)}$ is the sum of the $k$ largest singular values of $A$. For three different instructive proofs we refer to [2], [23] and [25, p.\ 56]. We also recall that
$$
\Vert A\Vert_{(k)}=\max_E \Vert AE\Vert_1
$$
where $E$ runs over the set of rank $k$ projections and $\Vert\cdot\Vert_1$ is the trace norm.

\vskip 10pt\noindent
{\bf Proof of Lemma 3.3.} Let $\{e_i\}$ be an orthonormal basis of eigenvectors for $A$ and $\{a_i\}$ the corresponding aeigenvalues. Letting $\Vert\cdot\Vert_p$ denote Schatten $p$-norms and using the fact that the norm of the  diagonal  is less than the norm of the full matrix, 
\begin{align*}
{\rm Tr\,}A^pZ^pA^p- {\rm Tr\,}(AZA)^p &= {\rm Tr\,}A^pZ^pA^p- \Vert AZA\Vert_p^p  \\
&\le\sum_i a_i^{2p}\langle e_i, Z^pe_i\rangle -\sum_i a_i^{2p}\langle e_i, Ze_i\rangle^p \\
&\le C(a,b,p) {\rm Tr\,}A^{2p}
\end{align*}
 where the second inequality follows from (18). \qquad $\Box$

\vskip 10pt\noindent
{\bf Proof of Theorem  3.4.} The main step consists in showing that the result holds for each Ky Fan $k$-norm,
\begin{equation}
\Vert A^pZ^pA^p\Vert_{(k)} - \Vert (AZA)^p\Vert_{(k)} \le C(a,b,p) \Vert A^{2p}\Vert_{(k)}.
\end{equation}
To this end, note that there exists a rank $k$ projection $E$ such that
\begin{align*}
\Vert A^pZ^pA^p\Vert_{(k)} &=\Vert EA^pZ^pA^pE\Vert_1 \\
&=\Vert Z^{p/2}A^pEA^pZ^{p/2}\Vert_1 \\
&=\Vert (A^pEA^p)^{1/2}Z^p(A^pEA^p)^{1/2}\Vert_1.
\end{align*}
We may then apply Lemma 3.3 to get
$$
\Vert A^pZ^pA^p\Vert_{(k)} \le \Vert \{ (A^pEA^p)^{1/2p}Z(A^pEA^p)^{1/2p} \}^p \Vert_1 +C(a,b,p) \Vert A^pEA^p\Vert_1.
$$
Since $\Vert A^pE^pA^p\Vert_1=\Vert EA^{2p}E\Vert_1 \le \Vert A^{2p}\Vert_{(k)}$ it suffices to show
$$
\Vert \{ (A^pEA^p)^{1/2p}Z(A^pEA^p)^{1/2p} \}^p \Vert_1 \le \Vert (AZA)^p\Vert_{(k)}
$$
or equivalently
\begin{equation}
\Vert \{ Z^{1/2}(A^pEA^p)^{1/p} Z^{1/2} \}^p \Vert_1 \le \Vert ( Z^{1/2}A^2 Z^{1/2})^p\Vert_{(k)}. 
\end{equation}
Let $X=Z^{1/2}(A^pEA^p)^{1/p} Z^{1/2}$ and $Y=Z^{1/2}A^2 Z^{1/2}$. Since $t\longrightarrow t^{1/p}$ is operator monotone
we infer $X\le Y$. Next we note that there exists a rank $k$ projection $F$ such that $FX=XF$ and $\Vert X^p\Vert_1=\Vert FX^pF\Vert_1$. Hence we may apply the auxillary lemma below to obtain
$$
\Vert X^p\Vert_1\le\Vert FY^pF\Vert_1
$$
which is the same as (20).

Having proved (19), let us show the general case. We first write (19) as
\begin{equation}
\Vert A^pZ^pA^p\Vert_{(k)} \le \Vert (AZA)^p\Vert_{(k)} + C(a,b,p) \Vert A^{2p}\Vert_{(k)}
\end{equation}
and we introduce a unitary $V$ such that $(AZA)^p$ and $VA^{2p}V^*$ form a monotone pair. Then (21) is equivalent to 
\begin{align*}
\Vert A^pZ^pA^p\Vert_{(k)} &\le \Vert (AZA)^p\Vert_{(k)} + C(a,b,p) \Vert VA^{2p}V^*\Vert_{(k)} \\
&=\Vert (AZA)^p + C(a,b,p)  VA^{2p}V^*\Vert_{(k)}
\end{align*}
by the simple fact that $\Vert X+Y\Vert_{(k)}=\Vert X\Vert_{(k)}+\Vert Y\Vert_{(k)}$ for all positive monotone pairs $(X,Y)$.
Therefore, Fan's dominance principle entails
$$
\Vert A^pZ^pA^p\Vert \le \Vert (AZA)^p + C(a,b,p) VA^{2p}V^*\Vert
$$
and the triangular inequality completes the proof. \qquad $\Box$

\vskip 10pt\noindent
{\bf Lemma 3.5.} {\it Let $0\le X\le Y$ and let $F$ be a projection, $FX=XF$. Then, for all $p>1$,
$$
{\rm Tr}\,FX^pF \le {\rm Tr}\, FY^pF.
$$
}

\noindent
{\bf Proof.} Compute ${\rm Tr}\,FX^pF$ in a basis of eigenvectors for $XF$ and apply (16). \qquad $\Box$

\vskip 10pt\noindent
{\bf Remark.} In [14, 16] several results related to Theorem 3.2, 3.4 are given for the operator norm. For instance,
 in [16] the authors prove Theorem 3.4 ([16, Corollary 9]) and give results for $0<p<1$. A special case of resuls in [14] is an additive version of Lemma 2.2. Under the same assumptions of this lemma, the authors show:
$$
\Vert AZ\Vert_{\infty}-\rho(AZ) \le \frac{(a-b)^2}{4(a+b)}\Vert A\Vert_{\infty}.
$$
Of course, letting $A$ be a rank one projection we recapture a classical reverse inequality, companion to (4).
Finally, let us mention a new book in which the reader may find many other reverse inequalities and references,  T.\ Furuta, J.\ Mi\'ci\'c, J.\ Pe$\check{{\rm c}}$ari\'c and Y. Seo, {\it Mond-Pe$\check{{\rm c}}$ari\'c Method in Operator Inequalities,}  Monograph in Inequalities 1, Element, Zagreb, 2005.

\vskip 20pt
{\bf References}
\vskip 2pt
\noindent
{\small 
\noindent
[1] H.\ Araki, {\it On an inequality of Lieb and Thirring,} Letters in Math.\ Phys. 19 (1990) 167-170.
\vskip 2pt
\noindent
[2] R.\ Bhatia, Matrix Analysis, Springer, Germany, 1996.
\vskip 2pt
\noindent
[3] J.-C.\ Bourin,  Some  inequalities for norms on matrices and operators, Linear Alg.\ Appl.\ 292 (1999) 139-154.
\vskip 2pt
\noindent
[4] J.-C.\ Bourin,  Singular values of compressions, dilations and restrictions, Linear Alg.\ Appl.\ 360 (2003) 259-272.

\vskip 2pt
\noindent
[5] J.-C.\ Bourin, {\it Symmetric  norms and reverse inequalities to Davis and Hansen-Pedersen characterizations of operator convexity,} to appear in Math.\ Ineq.\ Appl.
\vskip 2pt
\noindent
[6] J.-C.\ Bourin, Compressions, Dilations and Matrix Inequalities, RGMIA monograph, Victoria university, 
Melbourne 2004 (http://rgmia.vu.edu.au/monograph).
\vskip 2pt
\noindent
[7] J.-C.\ Bourin,  Reverse inequality to Araki's inequality: Comparison of $A^pZ^pZ^p$ and $(AZA)^p$  Math.\ Ineq.\ Appl. 8${\rm n^0}$3 (2005) 373-378.
\vskip 2pt
\noindent
[8] J.-C.\ Bourin,  Singular values of products of positive operators $AZB$ and $ZAB$, Linear Alg.\ Appl.\ 407 (2005) 64-70.
\vskip 2pt
\noindent
[9] J.-C.\ Bourin,  Reverse rearrangement inequalities via matrix technics, J.\ Inequal. Pure and Appl., submitted.

\vskip 2pt\noindent
 [10] S.S. Dragomir, {\it Reverse of Shwarz, triangle and Bessel Inequalities}, RGMIA Res.\ Rep.\ Coll.\ 6 (2003), Suppl.\ article 19.
 
\vskip 2pt\noindent
 [11] S.S. Dragomir, {\it A survey on Cauchy-Bunyakowsky-Schwarz type discrete inequalities}, J.\ Inequal.\ Pure and Appl.\ Math.\ 4$n^0$3 Art.\ 63, 2003.
 
 \vskip 2pt\noindent
 [12] N.\ Elezovic, L.\ Marangunic, J.\ Pecaric,  Unified treatement of complemented Schwarz and Gruss inequalities in inner product spaces, Math.\ Ineq.\ appl., 8${\rm n^0}$2 (2005), 223-231
 
\vskip 2pt\noindent
\noindent
[13] J.I.\ Fujii, A trace inequality arising from quantum information theory, Linear Alg.\ Appl.\ 400 (2005) 141-146.
 
\vskip 2pt
\noindent
[14] J.I.\ Fujii, Y.\ Seo and M.\ Tominaga, {\it Kantorovich type inequalities for operator norm,} Math.\ Ineq.\ Appl.,
8${\rm n^0}$3 (2005) 529-535.
\vskip 2pt
\noindent
[15] M.\ Fujii, T.\ Furuta, R.\ Nakamoto and S.\ E.\ Takahasi, Operator inequalities  and covariance in noncommutative probability space, Math.\ Japon.\ 46 (1997) 317-320.
\vskip 2pt
\noindent
[16] M.\ Fujii and Y.\ Seo,  Reverse inequalities of Araki, Cordes and Lowner-Heinz inequalities, to appear in Math.\ Ineq.\ Appl.
\vskip 2pt
\noindent
[17] T.\, Furuta,  Specht ratio S(1) can be expressed by Kantorovich constant $K(p)$: $S(1)=exp[K'(1)]$ and its applications,  Math.\ Ineq.\ Appl. 6 ${\rm n}^0$3 (2003) 521-530.
\vskip 2pt
\noindent
[18] T.\ Furuta and  Y.\ Seo,   An application of generalized Furuta inequality to Kantorovich type inequalities, Sci.\ Math.\ 2 (1999) 393-399.
\vskip 2pt
\noindent
[19] S.\ Izumino, H.\ Mori and Y.\ Seo, On Ozeki's inequality, J.\  Inequal.\  Appl., 2 (1998) 235-253.
\vskip 2pt
\noindent
[20] Ky Fan,  Some matrix inequalities, Abh.\ Math.\ Sem.\ Univ.\ Hamburg 29 (1966) 185-196.
\vskip 2pt
\noindent
[21] D.\ S.\ Mitrinovic, J.\ E.\ Pe$\check{{\rm c}}$ari\'c and A.\ M.\ Fink, Classical and new inequalities in Analysis, Mathematics and its applications (East european series) 61, Kluwer academic publisher, Dordrecht, 1993.
\vskip 2pt
\noindent
[22] B.\ Mond,  J.E.\ Pe$\check{{\rm c}}$ari\'c, A matrix version of the Ky Fan generalization of the Kantorovich inequality,
Linear and Multilinear Algebra 36 (1994) 217-221.
\vskip 2pt
\noindent
[23] B.\ Simon,  Trace ideals and their applications, LMS lecture note, 35 Cambridge Univ.\ Press, Cambridge, 1979.

\vskip 2pt
\noindent
[24] T.\ Yamazaki,  an extension of Specht's theorem via Kantorovich's inequality and related results,  Math.\ Ineq.\ Appl. 3 ${\rm n}^0$1 (2000) 89-96.
\vskip 2pt
\noindent
[25] X.\ Zhan,  Matrix inequalities, LNM 1790, Springer, 2002, Germany.

\end{document}